\documentclass[a4paper]{scrartcl}

\usepackage{scrpage2}
\usepackage{amsmath}
\usepackage{amsfonts}
\usepackage{epsfig}
\usepackage{graphics} 
\usepackage{url}

\date{}
\ihead{\small{\textsc{Ohlberger, Rave, Schindler}}}
\ohead{\small{\textsc{MOR for Multiscale Lithium-Ion Batteries}}}

\begin{document}

\title{Model Reduction for Multiscale Lithium-Ion Battery Simulation}

\author{Mario Ohlberger$^*$ \and Stephan Rave\thanks{Applied Mathematics M\"unster, CMTC \& Center for Nonlinear Science, University of M\"unster, Einsteinstr.\ 62, 48149 M\"unster, Germany, \texttt{mario.ohlberger,stephane.rave,felix.schindler} \texttt{@uni-muenster.de}. This work has been supported by the German Federal Ministry of Education and Research (BMBF) under contract number 05M13PMA.} \and Felix Schindler$^*$}

\maketitle

\begin{abstract}
In this contribution we are concerned with efficient model reduction for multiscale problems 
arising in lithium-ion battery modeling with spatially resolved porous electrodes. We present new
results on the application of the reduced basis method to the resulting instationary 3D battery model 
that involves strong non-linearities due to Buttler-Volmer kinetics. Empirical operator interpolation is used 
to efficiently deal with this issue. Furthermore, we present the localized reduced basis multiscale method 
for parabolic problems applied to a thermal model of batteries with resolved porous electrodes. 
Numerical experiments are given that demonstrate the reduction capabilities of the presented approaches
for these real world applications.
\end{abstract}

\section{Introduction}
\label{sec:1}

Continuum modeling of batteries results in a reaction-diffusion-transport system of coupled nonlinear partial
differential equations in complex multiscale and multi-phase pore structures. In recent contributions 
\cite{LatzZausch2011,PopovVutovEtAl2011,LessSeoEtAl2012}
three dimensional numerical models have been proposed that resolve the porous electrodes and thus serve as a basis 
for multiscale modeling as well as for more complex modeling of degradation processes such as Lithium plating.
Concerning multiscale modeling in the context of battery simulation, we refer e.g. to \cite{MR2832465,MR2990422,Taralova2015}.
These models result in huge time dependent discrete systems which require enormous computing resources, already
for single simulation runs. Parameter studies, design optimization or optimal control, however, require
many forward simulation runs with varying material or state parameters and are thus virtually impossible.
Hence, model reduction approaches for the resulting parameterized systems are indispensable for such simulation 
tasks. In this contribution we apply the reduced basis method, that has seen significant advance in recent years.
For an overview, we refer to the recent monographs~\cite{HesthavenRozzaEtAl2016,QuarteroniManzoniEtAl2016}
and the tutorial~\cite{Ha14}.

Concerning model reduction for lithium-ion battery models, we refer to the early work \cite{CaiWhi:2009}  where 
Galerkin projection into a subspace generated by proper orthogonal decomposition (POD)  is used on the basis 
of the mathematical model proposed in \cite{Doyle19931526}. In \cite{LV15}, the POD approach is used in the 
context of parameter identification for battery models. Preliminary results concerning model reduction with 
reduced basis methods can be found in \cite{IlievLatzEtAl2012,VW13} and \cite{ORSZ14}.

In this contribution we focus on two advances in reduced order modeling for batteries. 
First, in Section \ref{sec:2}, we present new results concerning nonlinear model reduction 
for the microscale battery model presented in \cite{LatzZausch2011}. 
The model reduction approach is based on Galerkin projection onto POD spaces,
extended to nonlinear problems using empirical operator 
interpolation \cite{BMNP04,MR2537223,DHO12}.

Second, in Section \ref{sec:3} we demonstrate the applicability of the localized reduced basis multiscale method 
(LRBMS) for a thermal model of batteries with resolved porous electrodes. 
The LRBMS has first been introduced in \cite{KOH2011,Albrecht2012lq} and further developed in 
\cite{OS14,OS15}. The later contributions in particular propose a rigorous a posteriori error estimate 
for the reduced solution with respect to the exact solution for elliptic problems that is localizable and 
can thus be used to steer an adaptive online enrichment procedure. For an application of the method for more complex problems
in the context of two phase flow in porous media we refer to \cite{Kaulmann2015wa}

\section{Reduced basis methods applied to pore-scale battery models}
\label{sec:2}

In this section we present first numerical results for the full model order
reduction of large 3D pore-scale Li-ion battery models.
These results extend our preliminary findings in \cite{ORSZ14}, where we
tested the quality of the reduced basis approximation for a small test geometry,
towards realistically sized geometries used in real-world
simulations, showing the feasibility of our model reduction approach.
Before discussing our new results, we will briefly review the battery model
under consideration and the basics of the reduced basis methodology.

\subsection{A pore-scale Lithium-Ion battery model}

\begin{figure}[t]
	\centering
	\includegraphics{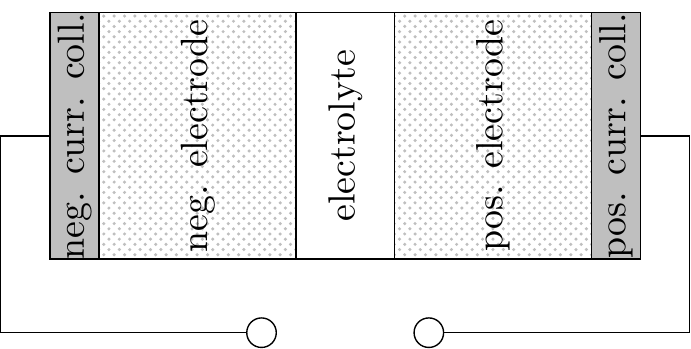}
	\caption{Schematic overview of the considered battery geometry (note that electrodes have
		porous structure, pore space is filled with electrolyte).}
	\label{fig:battery_schematic}
\end{figure}

Following \cite{ORSZ14}, we consider a pore-scale battery model based on
\cite{LatzZausch2011}.
The computational
domain is divided into five parts: electrolyte, positive/negative electrode,
positive/negative current collector (Fig.~\ref{fig:battery_schematic}).
On each of these subdomains, partial
differential equations are given for the Li-ion
concentration $c$ and the electrical potential $\phi$.

For the electrolyte we have
\begin{align}
\frac{\partial c}{\partial t} - \nabla \cdot (D_e \nabla c) &= 0,
\label{eq:c_electrolyte} \\
 -\nabla \cdot \Bigl(\kappa \frac{1 - t_+}{F}RT \frac{1}{c} \nabla c - \kappa \nabla
 \phi\Bigr) &= 0,
\end{align}
where $D_e = 1.622\cdot10^{-6}\frac{cm^2}{s}$, $\kappa = 0.02 \frac{s}{cm}$,
$t_+ = 0.39989$ denote the collective interdiffusion coefficient
in the electrolyte, the ion conductivity, and the transference number.
$R = 8.314 \frac{J}{mol\,K}$, $F = 96487 \frac{As}{mol}$ are the universal gas constant and the Faraday
constant.
We fix the global temperature $T$ to $298K$.

In the electrodes, $c$ and $\phi$ satisfy
\begin{align}
\frac{\partial c}{\partial t} - \nabla \cdot (D_s \nabla c) &= 0,
\label{eq:c_electrodes} \\
-\nabla \cdot (\sigma \nabla \phi) &= 0,
\end{align}
where $D_s = 10^{-10}\frac{cm^2}{s}$ is the ion diffusion coefficient in the electrodes,
and $\sigma = 10 \frac{s}{cm}$ ($\sigma = 0.38 \frac{s}{cm}$) in the negative
(positive) electrode denotes the electronic conductivity.

Finally, no Li-ions can enter the current collectors, so $c = 0$ on the whole current
collector subdomains. Moreover, $\phi$ again satisfies
\begin{equation}
	-\nabla \cdot (\sigma \nabla \phi) = 0,
\end{equation}
with $\sigma = 10 \frac{s}{cm}$ ($\sigma = 0.38 \frac{s}{cm}$) for the negative
(positive) current collector.

Note that for this in comparison to \cite{LatzZausch2011} slightly simplified
model (assuming constant $t_+$), the equations (\ref{eq:c_electrolyte}),
(\ref{eq:c_electrodes}) are linear and decoupled from the potential equations.
However, the coupling between the two variables is established by the interface
conditions at the electrode-electrolyte interfaces,
where the so-called Butler-Volmer kinetics are assumed: the
electric current (ion flux) $j$ (N) from the electrodes into the electrolyte is 
given by
\begin{equation}
       j = 2k\sqrt{c_ec_s(c_{max}- c_s)} \sinh\left(\frac{\phi_s - \phi_e -
		       U_0(\frac{c_s}{c_{max}})}{2RT} \cdot F \right),
       \quad
       N = \frac{j}{F}.
\end{equation}
Here, $c_{e/s}$ ($\phi_{e/s}$) denotes the Li-ion concentration (electrical
potential) at the electrolyte/electrode side of the interface.
$c_{max} = 24681\cdot 10^{-6}\frac{mol}{cm^3}$ ($c_{max} = 23671 \cdot
10^{-6}\frac{mol}{cm^3}$) denotes the maximum Li-ion concentration in the
negative (positive) electrode, and the rate constant $k$ is given by $k = 0.002
\frac{A cm^{2.5}}{mol^{1.5}}$
at the negative and by $k = 0.2 \frac{A cm^{2.5}}{mol^{1.5}}$ at the positive electrode interface.
Finally, the open circuit potential is given by $U_0(s) = (-0.132+1.41\cdot
e^{-3.52 s})V$ for the negative, and by
\begin{equation}
\begin{aligned}
	    U_0(s) = &\Bigl[\ 0.0677504 \cdot \tanh (-21.8502\cdot s+ 12.8268) \\
	             &\ - 0.105734 \cdot \bigl( (1.00167 - s)^{-0.379571} - 1.576\bigr) \\
	             &\ - 0.045\cdot e^{-71.69 \cdot s^8} + 0.01\cdot e^{-200\cdot(s-0.19)}
		     + 4.06279\ \Bigr]\cdot V
\end{aligned}
\end{equation}
for the positive electrode.

Given the porous electrode structures, these interface conditions apply to a
large surface area, giving this model highly nonlinear
dynamics.

Finally, the system is closed by the following boundary conditions:
homogeneous Neumann conditions for $c$ at all further inner
and external domain boundaries, continuity conditions for $\phi$ at the current
collector-electrode interfaces, homogenous Neumann conditions for $\phi$ at
the current collector-electrolyte interfaces, $\phi \equiv U_0(c(0) / c_{max})$ at
the negative current collector boundary, and $-n\cdot \sigma \nabla \phi \equiv
\mu$ at the positive current collector boundary.

We consider the fixed charge rate $\mu$ as a parameter we want to vary in our
numerical experiments.

\subsection{Reduced basis method and empirical interpolation}

After cell-centered finite volume discretization of the model on a voxel grid,
replacing the numerical fluxes by the Butler-Volmer relations at the
electrode-electrolyte interfaces, and backward Euler time discretization,
we arrive at nonlinear, discrete equations systems of the form
\begin{equation}
\label{eq:detailed}
        \begin{bmatrix}
                 \frac{1}{\Delta t}(c_{\mu}^{(t+1)} - c_{\mu}^{(t)}) \\
                 0
        \end{bmatrix}
         + A_\mu
         \left(\begin{bmatrix}
                c_{\mu}^{(t+1)} \\
                \phi_{\mu}^{(t+1)}
         \end{bmatrix}\right)
         = 0, \qquad (c_{\mu}^{(t)}, \phi_{\mu}^{(t)}) \in V_h \oplus V_h,
\end{equation}
where $A_\mu$ denotes the parametric finite volume space differential operator
acting on the finite volume space $V_h$ (see \cite{PopovVutovEtAl2011} for a detailed derivation).
Solving these systems using Newton's method requires many hours for realistic
geometries, even when using advanced algebraic multigrid solvers for computing
the Newton updates.

Projection-based parametric model reduction methods are based on the idea of
finding problem adapted approximation spaces $\tilde{V} \subseteq V_h \oplus V_h$
in which a reduced order solution is obtained by projection of the original
equation system:
\begin{equation}
\label{eq:reduced}
		P_{\tilde{V}} \left\{
                \begin{bmatrix}
                         \frac{1}{\Delta t}(\tilde{c}_{\mu}^{(t+1)} - \tilde{c}_{\mu}^{(t)}) \\
                         0
                \end{bmatrix} 
                 + 
                 A_\mu
                 \left(\begin{bmatrix}
		        \tilde{c}_{\mu}^{(t+1)} \\
			\tilde{\phi}_{\mu}^{(t+1)}
                 \end{bmatrix}\right) \right\}
         = 0, \quad (\tilde{c}_{\mu}^{(t)}, \tilde{\phi}_{\mu}^{(t)}) \in
	 \tilde{V}.
\end{equation}
Here, $P_{\tilde{V}}$ denotes the orthogonal projection onto $\tilde{V}$.
Since the manifold of system states $\{(c_\mu^{(t)}, \phi_\mu^{(t)})\ |\ \mu \in [\mu_{min},
\mu_{max}],\ t \in \{0, \ldots, T\} \}$ has a low-dimensional para\-me\-tri\-zation (by $(\mu, t) \in
\mathbb{R}^2)$, and assuming that this para\-metrization is sufficiently smooth,
there is hope to find low-dimensional approximation spaces
$\tilde{V}$ such that the model reduction error between the reduced solutions
(\ref{eq:reduced}) and the corresponding high-dimensional solutions
(\ref{eq:detailed}) is very small.

A vast amount of methods for constructing reduced spaces $\tilde{V}$ has been
considered in literature.
For time-dependent problems, the \textsc{POD-Greedy} method 
\cite{MR2405149,Ha13} has shown to produce approximation spaces with
quasi-optimal $l^\infty$-in-$\mu$, $l^2$-in-time reduction error.
In our experiments below, we apply a more basic approach by computing a basis
for $\tilde{V}$ via PODs of a pre-selected set of solution trajectories of
(\ref{eq:detailed}).
More precisely, we compute separate reduced concentration ($\tilde{V}_c$) and potential
($\tilde{V}_\phi$) spaces and let $\tilde{V}:= \tilde{V}_c \oplus
\tilde{V}_\phi$.
Due to the basic properties of POD, $\tilde{V}_c$, $\tilde{V}_\phi$ are
$l^2$-in-$\mu$, $l^2$-in-time best-approximation spaces for the considered 
training set of solutions.

Even though the equation systems (\ref{eq:reduced}) are posed on the
low-dimensional space $\tilde{V}$, solving (\ref{eq:reduced}) requires
the evaluation of the projected operator $P_{\tilde{V}} \circ A_\mu$ (and its
Jacobian), which in turn makes the computationally expensive evaluation of
$A_\mu$ on the full finite volume space $V_h \oplus V_h$ necessary.
The method of choice to overcome this limitation for nonlinear operators $A_\mu$
is empirical operator interpolation:
$A_\mu$ is replaced by an interpolant
$I_M \circ \tilde{A}_{M, \mu} \circ R_{M^\prime}$, where
$\tilde{A}_{M, \mu}: \mathbb{R}^{M^\prime} \to \mathbb{R}^M$ is the restriction
of $A_\mu$ to $M$ appropriately selected degrees of freedom (DOFs),
$R_{M^\prime}: V_h \oplus V_h \to \mathbb{R}^{M\prime}$ is the restriction of
the finite volume vectors to the $M^\prime$ DOFs required for the evaluation of
$\tilde{A}_{M, \mu}$ and $I_M: \mathbb{R}^M \to V_h \oplus V_h$ is the linear
combination with an appropriate interpolation basis (collateral basis).
Due to the locality of finite volume operators, $M^\prime$ can be chosen such
that $M^\prime \leq C\cdot M$, where $C$ only depends on the maximum number of
neighboring cells in the given mesh.
The interpolation DOFs and the associated collateral basis are obtained
from solution snapshot data using the \textsc{EI-Greedy} algorithm \cite{MR2537223,DHO12}.

A direct application of this approach to $A_\mu$ would not be
successful, however: since the collateral basis is contained in the linear span of
operator evaluations on solution trajectories, the $\phi$-parts of the collateral
basis vectors would, according to (\ref{eq:detailed}), completely vanish.
Therefore, we first decompose $A_\mu$ as $A_\mu = A^{(const)} + \mu \cdot A^{(bnd)} +
A^{(lin)} + A^{(1/c)} + A^{(bv)}$, where $A^{(1/c)}$, $A^{(bv)}$ are the
nonlinear operators corresponding to $-\nabla \cdot \kappa \frac{1 - t_+}{F}RT \frac{1}{c} \nabla c$
and the Butler-Volmer interfaces, $A^{(const)}$ ($A^{(bnd)}$) is the constant
(parametric) part of $A_\mu$ corresponding to the boundary conditions, and
$A^{(lin)}$ is the remaining linear part of $A_\mu$. We then apply empirical
operator interpolation separately to $A^{(1/c)}$ and $A^{(bv)}$. 
With $T[\tilde{c}^{(t)}_\mu](\tilde{c}, \tilde{\phi}):= (1/\Delta t\cdot (\tilde{c} - \tilde{c}_\mu^{(t)}),\ 0)$,
we arrive at the fully reduced systems
\begin{equation}\label{eq:ei_reduced}
\begin{aligned}
	\Bigl\{T[\tilde{c}^{(t)}_\mu] & +
		 P_{\tilde{V}} \circ A^{(const)} + \mu \cdot P_{\tilde{V}}
		 \circ A^{(bnd)} 
		  + P_{\tilde{V}} \circ A^{(lin)}  \\
		 & + \{P_{\tilde{V}} \circ I_{M^{(1/c)}}^{(1/c)}\} \circ
		 \tilde{A}_{M^{(1/c)}, \mu}^{(1/c)} \circ R_{M^{\prime (1/c)}}^{(1/c)} \\
		 & + \{P_{\tilde{V}} \circ I_{M^{(bv)}}^{(bv)}\} \circ
		 \tilde{A}_{M^{(bv)}, \mu}^{(bv)} \circ R_{M^{\prime (bv)}}^{(bv)}\  \ \  \Bigr\}
                 \left(\begin{bmatrix}
		        \tilde{c}_{\mu}^{(t+1)} \\
			\tilde{\phi}_{\mu}^{(t+1)}
                 \end{bmatrix}\right)
         = 0. 
 \end{aligned}
\end{equation}
After pre-computation of the linear maps $P_{\tilde{V}} \circ A^{(bnd)}$,
$P_{\tilde{V}} \circ A^{(lin)}$, $P_{\tilde{V}} \circ I_{M^{(1/c)}}^{(1/c)}$,
$R_{M^{\prime (1/c)}}^{(1/c)}$, $P_{\tilde{V}} \circ I_{M^{(bv)}}^{(bv)}$,
$R_{M^{\prime (bv)}}^{(bv)}$ and of the constant map $P_{\tilde{V}} \circ
A^{(const)}$ w.r.t.\ to a basis of $\tilde{V}$, (\ref{eq:ei_reduced}) 
can be solved quickly and independent of the dimension of $V_h$.

\subsection{Numerical experiments}
\label{sec:2.3}

\begin{figure}[t]
	\begin{center}
		\hspace{0.1cm}
		\raisebox{0.3\height}{\includegraphics[trim=0 80 80 20, clip,
			width=0.45\textwidth]{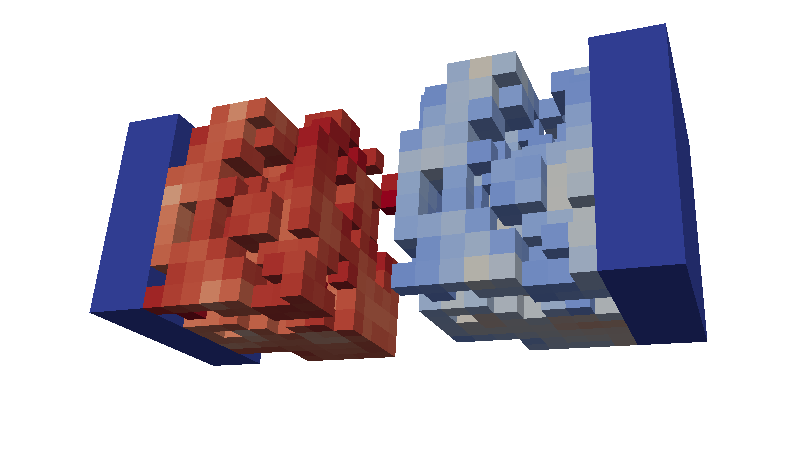}}
		\hfill
		\includegraphics{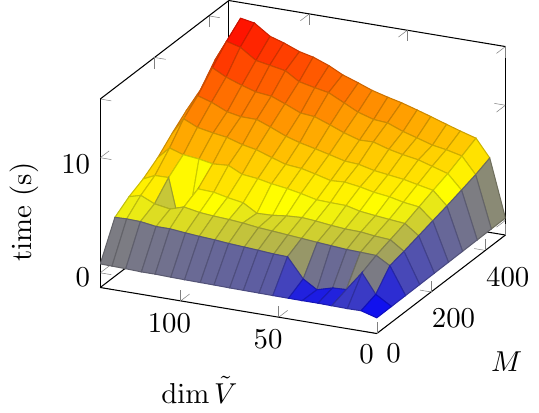}\\[\medskipamount]
    \includegraphics{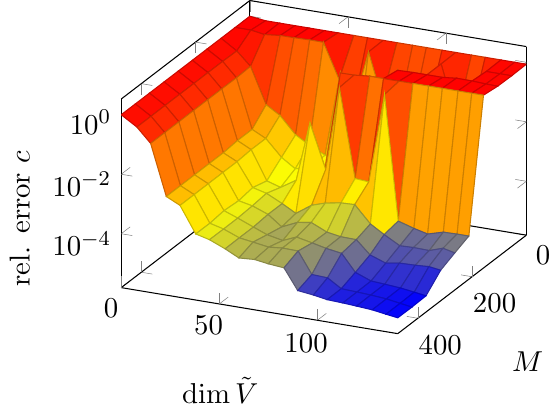}
		\hfill
		\includegraphics{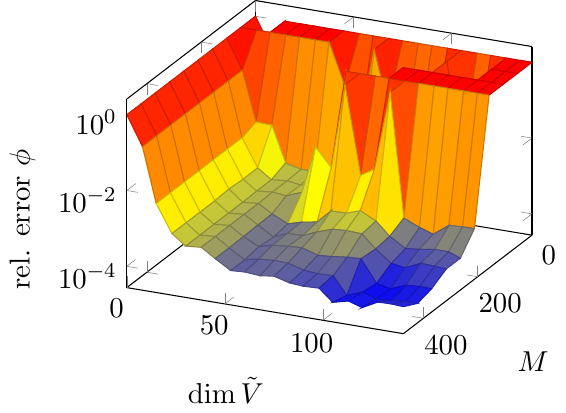}
	\end{center}
        \caption{%
Top left: small porous battery geometry used in numerical experiments.
Size: $104\mu m \times 40\mu m \times 40\mu m$, $4.600$ DOFs, coloring indicates Li$^+$
concentration at end of simulation, electrolyte not depicted.
Top right: average solution time in seconds vs.\ dimension of reduced space $\tilde{V}$
and number of interpolation points ($M:=M^{(1/c)}+ M^{(bv)}$).
Bottom: relative model reduction errors (\ref{eq:rel_error}) for concentration (left) and potential (right)
variable vs.\ dimension of reduced space and number of interpolation points.
A training set of 20 equidistant parameters was used for the generation of
$\tilde{V}$ and the interpolation data, $\# \mathcal{S}_{test} = 20$.
\label{fig:battery_small}
}
\end{figure}

\begin{figure}[t]
	\begin{center}
	\includegraphics[trim=80 120 60 80, clip, width=0.98\textwidth]{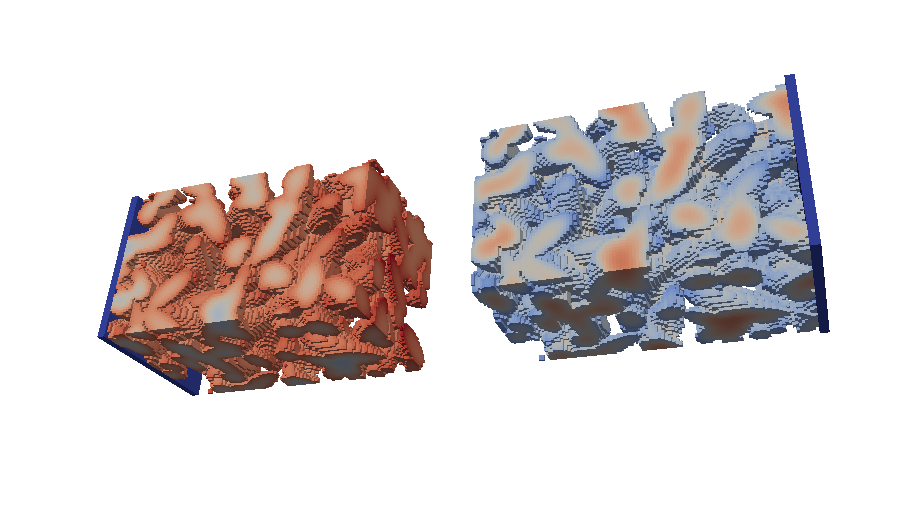}
	\end{center}
        \caption{
                Porous battery geometry used in the numerical experiments. Size:
		$246\mu m \times 60\mu m \times 60\mu m$, $1.749.600$ DOFs, coloring indicates
		Li$^+$ concentration at end of simulation, electrolyte not depicted.	\label{fig:simulation}}
\end{figure}

We consider two different test cases:
a small test geometry (Fig.~\ref{fig:battery_small}) which still exhibits the
most important properties of a real battery geometry, and a large, fully
resolved geometry (Fig.~\ref{fig:simulation}) useable for real-world
simulations.
In both cases, the initial Li$^+$ concentration $c_0$ was set to $c_0 \equiv
2639 \cdot 10^{-6} \frac{mol}{cm^3}$ ($c_0 \equiv 20574 \cdot 10^{-6}
\frac{mol}{cm^3}$) for the positive (negative) electrode and to $c_0 \equiv 1200
\cdot 10^{-6} \frac{mol}{cm^3}$ in the electrolyte.
The model was simulated on a $T = 2000 s$ ($T = 1600 s$) time interval for the small
(large) geometry, with a time step size of $\Delta t = 20 s$.
The charge rate $\mu$ was for each simulation chosen as a constant from the
interval $\bigl[0.00012 \frac{A}{cm^2}, 0.0012 \frac{A}{cm^2}\bigr]$ for the small and
from the interval $\bigl[0.000318 \frac{A}{cm^2}, 0.00318 \frac{A}{cm^2}\bigr]$ for the
large geometry.

To generate the reduced space $\tilde{V}$, we computed solution snapshots on
training sets $\mathcal{S}_{train}$ of equidistant parameters.
For the small geometry we chose $\# \mathcal{S}_{train} = 20$, whereas for the
large geometry we only selected the lower and upper boundary of the considered
parameter domain, i.e.\ $\# \mathcal{S}_{train} = 2$.
For the generation of the empirical interpolation data using the
$\textsc{EI-Greedy}$ algorithm, we additionally included the evaluations of
$A_\mu^{(1/c)}$ and $A_\mu^{(bv)}$ on all intermediate Newton stages of the
selected solution trajectories.

As a measure for the model reduction error we consider the relative
$l^\infty$-in-$\mu$, $l^\infty$-in-time error given by
\begin{equation}\label{eq:rel_error}
\max_{\mu \in \mathcal{S}_{test}}\ \max_{t \in \{0, 1, \ldots
	T/\Delta t\}}\ \frac{\|u_\mu^{(t)} - \tilde{u}_\mu^{(t)} \|}{\max_{t \in \{0, 1, \ldots
	T/\Delta t\}} \|u_\mu^{(t)}\|},
\end{equation}
where $u$ ($\tilde{u}$) is the concentration or potential part of the (reduced) solution
and $\mathcal{S}_{test}$ denotes a random set of test parameters.

All simulations of the high-dimensional model have been performed with the
battery simulation software \texttt{BEST} \cite{LessSeoEtAl2012}, which has been
integrated with our model order reduction library \texttt{pyMOR} \cite{ORSZ14,MRS15}.
The experiments were conducted as single-threaded processes on a dual socket
compute server equipped with two Intel Xeon E5-2698 v3 CPUs with 16 cores
running at 2.30 GHz each and 256GB of memory available.

For the small test geometry, we observe a rapid decay of the model reduction
error for both the concentration and the potential variable (Fig.~\ref{fig:battery_small}).
As usual for empirical operator interpolation, we see that the number of
interpolation points has to be increased for larger reduced space dimensions
in order to ensure stability of the reduced model.
Doing so, we obtain relative reduction errors as small as $10^{-4}$ with
simulation times of less than $15s$.

Since we only selected 2 solution trajectories for the generation of the reduced
model for the large geometry, we cannot expect such small model reduction errors over
the whole parameter domain.
In fact, the error stagnates already for relatively small reduced space
dimensions (Table~\ref{tab:battery_results}).
Nevertheless, we easily achieve errors of less than one percent for a simulation
time of $80s$.
With an average solution time for the high-dimensional model of over 6 hours, we
achieve at this error a speedup factor of 285.

Note that the solution time of the reduced model is still
significantly larger than for the small geometry.
This can be attributed to the fact that the localized evaluation of $A_\mu^{(1/c)}$,
$A_\mu^{(bv)}$ has been only partially implemented in \texttt{BEST} and still
requires operations on high-dimensional data structures. 
After the implementation of localized operator evaluation in \texttt{BEST} has been finalized,
we expect even shorter simulation times.

\newcommand{\BatteryDetailedTime}{22979}
\newcommand{\BatteryCBDim}{188}

\begin{table}[t]
\centering
\caption{Relative model reduction errors (\ref{eq:rel_error}) and reduced simulation times
	for the large battery geometry (Fig.\ \ref{fig:simulation}).
	\BatteryCBDim\ interpolation points, average time for solution
	of the high-dimensional model: $\BatteryDetailedTime s$, $\#
	\mathcal{S}_{test} = 10$.}
\label{tab:battery_results}
\begin{tabular}{p{2cm}p{2cm}p{2cm}p{2cm}p{2cm}p{2cm}}
\hline\noalign{\smallskip}
$\dim \tilde{V}$ & 11 & 21 & 30 & 40 \\
\hline\noalign{\smallskip}
rel. error $c$ & $9.26\cdot 10^{-3}$ & $3.96\cdot 10^{-3}$ & $3.05\cdot 10^{-3}$ & $2.93\cdot 10^{-3}$ \\
rel. error $\phi$ & $2.07\cdot 10^{-3}$ & $1.50\cdot 10^{-3}$ & $1.46\cdot 10^{-3}$ & $1.26\cdot 10^{-3}$ \\
time (s) & 82 & 81 & 79 & 81 \\
speedup & 279 & 285 & 290 & 283 \\

\noalign{\smallskip}\hline\noalign{\smallskip}
\end{tabular}
\end{table}

\section{Localized reduced basis multiscale approximation of heat conduction}
\label{sec:3}

The microscale battery model in Section \ref{sec:2} is considered under the 
assumption of constant global temperature $T$.
In general, it is desirable to couple this model with a spatially resolved
model for the temperature distributions within the battery.
For the model reduction of such heat conduction in porous electrodes we present
a first application of the localized reduced basis multiscale Method (LRBMS) for parabolic PDEs.

In this first step we consider the simulation and model reduction of heat
conduction separately from what is presented in Section \ref{sec:2} as a basis
for a coupled simulation and model reduction in future work.

For an introduction of the LRBMS for elliptic parameterized multiscale problems 
and recent results concerning localized a posteriori error estimation 
and online enrichment, we refer to \cite{OS15}.

\subsection{A battery - heat conduction model with resolved electrode geometry}

We consider here the same spatially resolved 3D pore-scale battery geometry 
(cf. Fig. \ref{fig:simulation}) as in Section
\ref{sec:2}, where the computational domain is composed of five materials which
are of interest for thermal modeling, that is: electrolyte, positive/negative
electrode and positive/negative current collectors, each with possibly
different thermal conductivities.

As a simplified model for heat conductivity within a battery with spatially resolved electrodes, 
we consider a parabolic PDE for the temperature $T$ of the form
\begin{align}
\label{eq:heat}
  \frac{\partial T}{\partial t} - \nabla\cdot\big( D\, \nabla T \big) = Q,
\end{align}
together with suitable initial and boundary conditions.
Here $D$ denotes the space-dependent conductivity tensor, which is material specific and thus 
takes different values in the current collectors, the porous electrodes, the separator, and the 
electrolyte. Hence, $D$ inherits the highly heterogeneous structure of the porous electrodes 
and thus has an intrinsic multiscale character. 
In general, $Q$ collects all heat generating sources, such as heat generation
due to electrochemical reaction, reversible heat and ohmic heat, each of which
may in turn depend on the Li-ion concentration and the electric potential and
thus vary in space and time. These sources arise in particular due to the electrochemical 
reaction at the interface between the electrodes and the electrolyte and it is thus desirable to consider
the full 3D pore-scale battery model in order to get an insight into possible variations of 
the temperature within the battery. 
We refer, e.g. to  \cite{CWW2005,CaiWhi:2009} for a more detailed derivation of 
an energy balance equation for Lithium-Ion batteries and corresponding simulation schemes. 

Depending on the study in question, any of the sources, the thermal conductivity
or the initial or boundary values may depend on a low-dimensional parameter
vector $\mu$.

\subsection{Localization of reduced basis methods - LRBMS}

As a first step towards a realistic model we allow for parametric thermal
conductivities and presume stationary sources and boundary values.
Thus, a (spatial) discretization of \eqref{eq:heat} by a suitable
discretization scheme (such as finite volumes or continuous or discontinuous
Galerkin (DG) finite elements) and a backward Euler
time-discretization yield a set of linear equations of the form,
\begin{align}
\label{eq:discrete_heat}
  \tfrac{1}{\Delta t} M_h \big(T^{(t + 1)} - T^{(t)}\big)
    + B_{h, \mu}\, T^{(t + 1)} = Q_h,
  &&T^{(t + 1)} \in V_h,
\end{align}
to be solved in each time step, where $M_h$ and $B_{h, \mu}$ denote
the discrete $L^2$-inner product and parametric space differential operators
induced by the spatial discretization, respectively, which act on the
corresponding high-dimensional discrete space $V_h$.
In addition, $Q_h$ denotes the discrete representation of the source and
boundary values.

To obtain a reduced order model for the discrete heat conduction model (\ref{eq:discrete_heat}), 
we proceed in an analog way, as described in Section \ref{sec:2} above, by a Galerkin projection 
onto a problem
adapted reduced approximation space $\tilde{V} \subset V_h$.
Once $\tilde{V}$ is given, we obtain the set of reduced equations for each time
step:
\begin{align}
\label{eq:reduced_heat}
  \tfrac{1}{\Delta t} \tilde{M} \big(\tilde{T}^{(t + 1)} -
\tilde{T}^{(t)}\big)
    + \tilde{B}_\mu\, \tilde{T}^{(t + 1)} = \tilde{Q},
  &&\tilde{T}^{(t + 1)} \in \tilde{V},
\end{align}
where $\tilde{M}$, $\tilde{B}_\mu$ and $\tilde{Q}$ denote the reduced operators
and functionals, respectively, acting on the low-dimensional reduced space
$\tilde{V}$.
Since all operators and functionals arising in \eqref{eq:reduced_heat} are
affinely decomposable with respect to the low-dimensional parameter vector $\mu$
(given for instance the thermal conductivity as in Section \ref{sec:3.3})
and linear with respect to $\tilde{V}$, we can precompute their
respective evaluations in a computationally expensive offline step, e.g., by
$\tilde{M} = \underline{P_{\tilde{V}}}^\bot
\underline{M_h} \; \underline{P_{\tilde{V}}}$, where $\underline{M_h}$ and
$\underline{P_{\tilde{V}}}$, respectively, denote the matrix representations of
$M_h$ and of the orthogonal projection $P_{\tilde{V}}: V_h \to \tilde{V}$ with
respect to the basis of $V_h$.
Online, for each new input parameter $\mu$, we can then quickly solve the
reduced low-dimensional problem \eqref{eq:reduced_heat} to obtain a
low-dimensional representation of the temperature $\tilde{T}$, which can be
post-processed to obtain the original temperature $T$, if required, or a
derived quantity of interest.

As mentioned above, the problem adapted reduced space $\tilde{V}$ can be
adaptively generated by an iterative \textsc{POD-Greedy} procedure
\cite{MR2405149}: in each step
of the greedy algorithm, given an error estimate on the model reduction error, a
full
high-dimensional solution trajectory for the hitherto worst-approximated
parameter is computed and the most dominant \textsc{POD} modes of the
projection error of this trajectory are added to the reduced basis spanning
$\tilde{V}$.

This procedure has been shown to produce quasi-optimal low-dimensional reduced
order models which successfully capture the dynamics of the original
high-dimensional model \cite{Ha13}.
However, in the context of multiscale phenomena or highly resolved geometries,
such as the porous structures within a Li-ion battery, the computational cost
required to generate the reduced model can become unbearably large, even given
modern computing hardware.

As a remedy, the localized reduced basis multiscale method has been introduced
for stationary elliptic multiscale problems \cite{KOH2011,Albrecht2012lq} to
lower the computational burden of traditional RB methods by generating
several local reduced bases associated with a partitioning of the computational
domain.
The local quantities associated with these individual subdomains can be
projected independently in parallel.
In \cite{OS14,OS15}, the LRBMS was extended to additionally account for the
discretization error and to allow for an adaptive enrichment of the local
reduced approximation spaces, which may even eliminate the need for global
solution snapshots at all.

In this contribution, we demonstrate a first application of the LRBMS to
parabolic multiscale problems, such as spatially resolved heat conduction in a
Lithium-Ion battery.
We therefore discretize \eqref{eq:heat} locally by a standard finite element or
discontinuous Galerkin scheme independently in each subdomain of a given
partitioning of the
computational domain and couple the arising local operators, products and
functionals along these subdomains by symmetric weighted interior penalty
discontinuous Galerkin fluxes (cf. \cite{OS15} and the references
therein).
We use the resulting discretization to compute global solution snapshots during
the greedy algorithm, as detailed above.
However, instead of a single reduced basis with global support, we
iteratively generate local reduced bases on each subdomain by localizing the
solution trajectories with respect to each subdomain and by carrying out local
\textsc{POD}s for further localized compression in a post-processing step.

The resulting reduced space is then given as the direct sum of the local
reduced approximation spaces spanned by these local reduced bases.
Accordingly, we obtain the reduced problem \eqref{eq:reduced_heat} by
local Galerkin projections of the local operators and functionals and coupling
operators associated with each subdomain and its neighbor, yielding sparse
reduced quantities.

\subsection{Numerical experiments}
\label{sec:3.3}

\begin{figure}[t]
  \footnotesize%
  \centering%
  \includegraphics{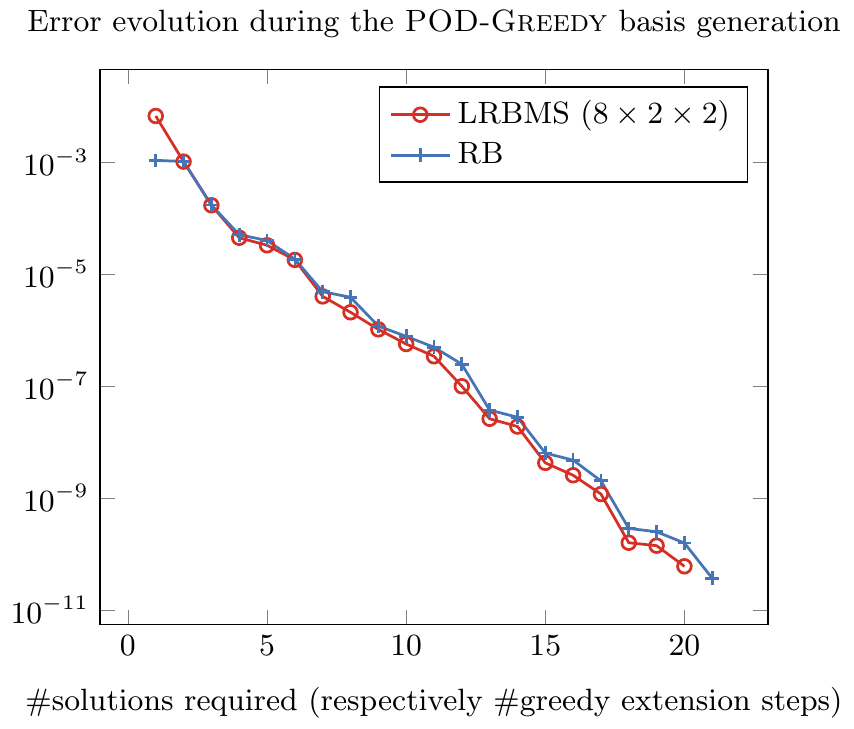}
  \caption{%
    \footnotesize%
Error evolution during the \textsc{POD-Greedy} basis generation to reach a
target absolute error of $10^{-10}$ for the numerical experiment from Section
\ref{sec:3.3}.
Depicted is the $L^\infty$-in-$\mu$, $L^\infty$-in-$t$, and $H^1$-in-space error
over the set of five equidistant training samples in $[0.1; 10]$.
}
  \label{fig:lrbms}
\end{figure}

To demonstrate the applicability of the LRBMS we conduct an experiment on the
same geometry used in the larger experiment in Section \ref{sec:2.3} (compare
Figure \ref{fig:simulation}).
For the thermal conductivities we choose constant values within each material
(the positive/negative electrode and the positive/negative current collectors),
as reported in \cite[4th column of Table 4]{CWW2005}.
Within the electrolyte we allow to vary the constant thermal conductivity
within the range $\mu \in [0.1; 10]$.
We pose homogeneous Dirichlet boundary values at the current collectors and
homogeneous Neumann boundary values elsewhere and start the simulations with an
initial temperature of 0K, using ten time steps to reach the final time
$10^{-3}$.
For the heat source we set $Q = 10^3$ within the electrodes and $Q = 0$ within
the current collectors and the electrolyte.
While this is not necessarily a physically meaningful setup, it inherits the
computational challenges of a realistic model, namely a highly resolved
geometry, discontinuous thermal conductivities depending on the materials and
heat sources which align with the geometry of the different materials.

We triangulate the computational domain with $5,313,600$ simplexes and compare
the LRBMS using $8 \times 2 \times 2$ subdomains to a standard RB method (which
corresponds to choosing one subdomain).
Within each subdomain, we use the same SWIPDG discretization as for the
coupling, thus yielding comparable discretizations with $21,254,400$
degrees of freedom in both approaches.
As an estimate on the model reduction error we use the true $L^\infty$-in-time, 
$H^1$-in-space error.

The discretization is implemented within the \texttt{DUNE} numerics environment
\cite{MR2421579,MR2421580}, centered around
\texttt{dune-gdt} \cite{wwwdunegdt}: the \texttt{dune-stuff} \cite{wwwdunestuff}
module provides classes for vectors, matrices and linear solvers (for instance the
\texttt{bicgstab.amg.ilu0} solver used in these experiments), \texttt{dune-gdt}
provides the discretization building blocks (such as discrete function spaces,
operators, products and functionals), and
\texttt{dune-hdd}\footnote{\url{https://github.com/pymor/dune-hdd}} provides
parametric discretizations compatible with \texttt{pyMOR} \cite{MRS15}.
Finally, \texttt{dune-pymor}\footnote{\url{https://github.com/pymor/dune-pymor}}
is used, as it provides the \texttt{Python}-bindings and wrappers to integrate the
\texttt{DUNE}-code with our model reduction framework \texttt{pyMOR}.
The experiments were conducted on the same compute server as described in Section \ref{sec:2.3}.

As we observe from Fig. \ref{fig:lrbms}, both the LRBMS and the standard RB
method show comparable exponential error decay.
In general, the quality of the reduced spaces generated by the LRBMS is
slightly better, while requiring less detailed solution snapshots to reach the
same target error.

\begin{table}[t]
\caption{%
Comparison of runtimes of the experiments from Section \ref{sec:3.3}.
Setup time includes grid generation, subdomain partitioning and assembly of
operators, products and functionals.
\textsc{POD-Greedy} time includes error estimation, generation of the reduced basis and the
reduced basis projection.
The average time to solve the detailed problem is $2h28m5s$.
}
\label{tab:lrbms}
\begin{tabular}{p{1.5cm}p{2cm}p{2.8cm}p{3cm}p{2cm}}
\hline\noalign{\smallskip}
        & setup    & \textsc{POD-Greedy} & reduced basis size & solution time\\
\noalign{\smallskip}\hline\noalign{\smallskip}
  RB    & $26m47s$ & $14h41m52s$         & 21                 & $35s$\\
  LRBMS & $36m7s$  & $14h34m39s$         & $32 \times 20$     & $35s$
\end{tabular}
\end{table}

As can be seen from Table \ref{tab:lrbms}, the \textsc{POD-Greedy} basis
generation using 32 subdomains is slightly faster than the basis generation
using a single subdomain.
However, since the experiments were conducted as single-threaded processes and
since the LRBMS allows for parallel local \textsc{POD}s and parallel local
reduced basis projections, the basis generation time of the LRBMS can be further
accelerated significantly.

\section{Conclusion}
\label{sec:conclusion}
In this contribution we have demonstrated the efficient applicability of recent model reduction approaches, 
such as the \textsc{POD-Greedy} reduced basis method, the empirical operator interpolation, and the 
localized reduced basis multiscale method (LRBMS) for efficient simulation of real world problems,
such as 3D spatially resolved heterogeneous Lithium-Ion battery models. 
The demonstrated model reduction approaches are realized within our model order reduction library 
\texttt{pyMOR} \cite{ORSZ14,MRS15} with bindings, both to the battery simulation software
 \texttt{BEST} \cite{LessSeoEtAl2012}, and the general purpose 
Distributed and Unified Numerics Environment \texttt{DUNE} \cite{MR2421579,MR2421580}, employing 
the \texttt{dune-gdt}, \texttt{dune-stuff}, and \texttt{dune-hdd} discretization and solver backends.
Speedup factors of about 285 were obtained for the full strongly non-linear battery model in Section \ref{sec:2}
using the reduced basis method with empirical operator interpolation \cite{DHO12}, and  around 253 for the linear parabolic
heat conduction model in Section \ref{sec:3} using a parabolic extension of the localized reduced basis multiscale method \cite{OS15}.

\section*{Acknowledgement}
The authors thank Sebastian Schmidt from Fraunhofer ITWM Kaiserslautern for the close and fruitful
collaboration within the BMBF-project MULTIBAT towards integration of \texttt{BEST} with \texttt{pyMOR}.

\ifx\undefined\bysame
\newcommand{\bysame}{\leavevmode\hbox to3em{\hrulefill}\,}
\fi

\end{document}